\documentclass[12pt]{elsart}
\usepackage{amsmath,amsfonts,amssymb,amsbsy}
\usepackage[cp1251]{inputenc}
\usepackage[english]{babel}
\usepackage{color}
\usepackage[pdftex]{graphicx}

\begin{document}

\begin{frontmatter}
\title{
Contact problem for a thin elastic layer with variable thickness:
Application to sensitivity analysis of articular contact mechanics
}
\author{I.~Argatov}
\ead{iva1@aber.ac.uk}
\corauth[cor]{Corresponding author.}
\address{Institute of Mathematics and Physics, Aberystwyth University,
Ceredigion SY23 3BZ, Wales, UK}
\begin{abstract}
In the framework of the recently developed asymptotic models for tibio-femoral contact incorporating  frictionless elliptical contact interaction between thin elastic, viscoelastic, or biphasic cartilage layers, we apply an asymptotic modeling approach for analytical evaluating the sensitivity of crucial parameters in joint contact mechanics due to small variations in the thicknesses of the contacting cartilage layers. The four term asymptotic expansion for the normal displacement at the contact surface is explicitly derived, which recovers the corresponding solution obtained previously for the 2D case in the compressible case. It was found that to minimize the influence of the cartilage thickness non-uniformity on the force-displacement relationship, the effective thicknesses of articular layers should be determined from a special optimization criterion. 
\end{abstract}

\end{frontmatter}

\setcounter{equation}{0}
\setcounter{figure}{0}
\section{Introduction}
\label{0vSection1}

Contact problems involving transmission of forces across biological joints are of considerable practical importance and a number of numerical models for articular contact are available \cite{Han_et_al2005,Wilson_et_al2005}. At the same time, the necessity of analytical models becomes an important issue in developing improved understanding of load distribution in the normal and pathological joints, which affects the mechanical aspects of osteoarthritis \cite{Ateshian1994,Wu1996}. Also, analytical modelling of the distributed internal forces generated by articular contact in tibio-femoral joints is required in multibody dynamic simulations of physical exercise of a human skeleton \cite{Eberhard1999,KlodowskiRantalainenMikkola2011}. As a rule, analytical models of articular contact assume rigid bones and represent cartilage as a thin elastic layer of constant thickness resisting to deformation like a Winkler foundation consisting of a series of discrete springs with constant length and stiffness \cite{BeiFregly2004}. However, a subject-specific approach to articular contact mechanics requires developing patient-specific models for accurate predictions. Recently, a sensitivity analysis of finite element models of hip cartilage mechanics with respect to varying degrees of simplified geometry was performed in \cite{Anderson2010}. 

Based on the asymptotic analysis of the frictionless contact problem for a thin elastic layer bonded to a rigid substrate in the thin-layer limit \cite{Argatov2005,Barber1990}, the following asymptotic model for contact interaction of two thin incompressible layers was established \cite{Argatov2011mb}:
\begin{equation}
-(E_1^{-1} h_1^3 +E_2^{-1} h_2^3)\Delta_y p({\bf y})=\delta_0-\varphi({\bf y}), \quad {\bf y}\in\omega,
\label{0v(1.1)}
\end{equation}
\begin{equation}
p({\bf y})=0,\quad \frac{\partial p}{\partial n}({\bf y})=0,\quad {\bf y}\in\Gamma.
\label{0v(1.2)}
\end{equation}
Here, $p({\bf y})$ is the contact pressure density, $h_\alpha$ and $E_\alpha$ are the thickness and elastic modulus of the layer material, respectively, $\alpha=1,2$, $\Delta_y=\partial^2/\partial y_1^2+\partial^2/\partial y_2^2$ is the Laplace differential operator, $\delta_0$ is the vertical approach of the rigid substrates, 
$\varphi({\bf y})$ is the gap function defined as the distance between the layer surfaces in the vertical direction, $\omega$ is the contact area, $\Gamma$ is the contour of $\omega$, $\partial/\partial n$ is the normal derivative. 

It was shown \cite{ArgatovMishuris2010e,Ateshian1994,BarryHolmes2001} that the problem (\ref{0v(1.1)}), (\ref{0v(1.2)}) describes the instantaneous response of thin biphasic layers to dynamic and impact loading. In \cite{ArgatovMishuris2011ve}, the elastic model (\ref{0v(1.1)}), (\ref{0v(1.2)}) was generalized for the general viscoelastic case.

With respect to articular contact, a special interest represents the case when the subchondral bones are shaped as an elliptic paraboloid
\begin{equation}
\varphi({\bf y})=(2R_1)^{-1}y_1^2+(2R_2)^{-1}y_2^2
\label{0v(1.3)}
\end{equation}
with positive curvature radii $R_1$ and $R_2$. 

In the case (\ref{0v(1.3)}), the exact solution to the problem (\ref{0v(1.1)}), (\ref{0v(1.2)}) has the following form \cite{ArgatovMishuris2010e,Barber1990}:
\begin{equation}
p({\bf y})=p_0\biggl(1-\frac{y_1^2}{a_1^2}-\frac{y_2^2}{a_2^2}\biggr)^2.
\label{0v(1.4)}
\end{equation}
Integration of the pressure distribution (\ref{0v(1.4)}) over the elliptical contact region $\omega$ with the semi-axes $a_1$ and $a_2$ results in the following force-displacement relationship \cite{Argatov2011mb}:
\begin{equation}
P=\frac{\pi m}{3}M_P(s)R_1 R_2\delta_0^3.
\label{0v(1.5)}
\end{equation}
Here, $M_P(s)$ is a dimensionless factor depending on the aspect ratio $s=a_2/a_1$, and the coefficient $m$ is given by
\begin{equation}
m=(E_1^{-1} h_1^3 +E_2^{-1} h_2^3)^{-1}.
\label{0v(1.6)}
\end{equation}

The asymptotic model (\ref{0v(1.1)})\,--\,(\ref{0v(1.3)}) assumes that the cartilage layers have constant thicknesses, whereas it is well known \cite{Akiyama2012} that articular cartilage has a variable thickness as well as the surface of subchondral bone deviates from the ellipsoid shape \cite{Siu1996}.
A sensitivity of the model (\ref{0v(1.1)}), (\ref{0v(1.2)}) with respect to small perturbations of the gap function (\ref{0v(1.3)}) was performed in \cite{ArgatovMishuris2011p}. In particular, it was shown \cite{Argatov2011mb} that the influence of the gap function variation on the force-displacement relationship will be negligible if the effective geometrical characteristics $R_1$ and $R_2$ are determined by a least square method. 

To our knowledge, in the literature there is only one study \cite{VorovichPeninin1971} where the 2D case of contact problem for a thin elastic strip of variable thickness was solved by an asymptotic method under the assumption that Poisson's ratio of the strip material is not very close to $0{.}5$. At the same time, many asymptotic solutions were derived for an elastic layer of constant thickness both in the axisymmetric \cite{Aleksandrov1968,Chadwick2002,Matthewson1981} and the non-axisymmetric \cite{Argatov2005,Barber1990,Hlavacek2008,Jaffar1989} cases.

In the present paper, a three-dimensional unilateral contact problem for a thin elastic layer of variable thickness bonded to a rigid substrate is considered. Two cases are studied separately: (a) Poisson's
ratio of the layer material is not very close to $0{.}5$; (b) the layer material is incompressible with Poisson's ratio of $0{.}5$. After developing a refined asymptotic model, we apply sensitivity analysis to determine how ``sensitive'' is the mathematical model (\ref{0v(1.1)}), (\ref{0v(1.2)}) to variations in the values of the layer thicknesses. To be more precise we consider the term ``sensitivity'' in a broad sense by allowing variable layer thicknesses, whereas the original model deals with scalar parameters $h_1$ and $h_2$. 

\section{Contact problem for a thin elastic layer with variable thickness}
\label{0vSection2}

We consider a homogeneous, isotropic, linearly elastic layer with a plane contact surface, $x_3=0$, and a variable thickness, $H(x_1,x_2)$, firmly attached to an uneven rigid surface
\begin{equation}
x_3=H(x_1,x_2).
\label{0v(2.1)}
\end{equation}

In the absence of body forces, equations governing small deformations of the layer are
\begin{equation}
\frac{\partial\sigma_{1j}}{\partial x_1}+\frac{\partial\sigma_{2j}}{\partial x_2}+
\frac{\partial\sigma_{3j}}{\partial x_3}=0,\quad j=1,2,3,
\label{0v(2.2)}
\end{equation}
\begin{equation}
\sigma_{ij}=\lambda\delta_{ij}(\varepsilon_{11}+\varepsilon_{22}+\varepsilon_{33})
+2\mu\varepsilon_{ij},
\label{0v(2.3)}
\end{equation}
\begin{equation}
\varepsilon_{ij}=\frac{1}{2}\biggl(\frac{\partial u_i}{\partial x_j}
+\frac{\partial u_j}{\partial x_i}\biggr),
\label{0v(2.4)}
\end{equation}
where $\sigma_{ij}$ is the Cauchy stress tensor, $\lambda$ and $\mu$ are Lam\'e parameters for the layer material, $\varepsilon_{ij}$ is the infinitesimal strain tensor, $u_j$ is the displacement component along the $x_j$-axis, $\delta_{ij}$ is Kronecker's delta. 

\begin{figure}[h!]
    \centering
\vbox{
    \includegraphics [scale=0.5]{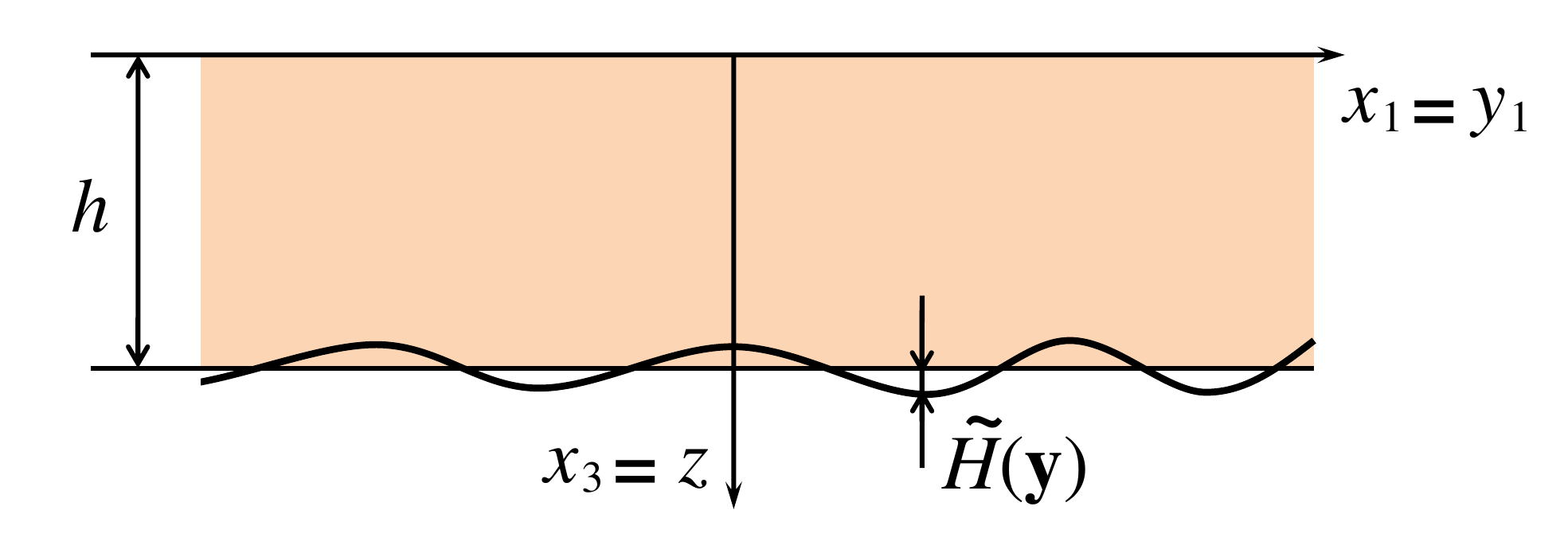}
    \caption{Elastic layer with a variable thickness.}
    \label{figure_strip.pdf}    
   }
\end{figure}

We assume that the elastic layer is indented by a smooth rigid punch in the form of an elliptic paraboloid
\begin{equation}
x_3=-\varphi(x_1,x_2), \quad \varphi(x_1,x_2)=(2R_1)^{-1}x_1^2+(2R_2)^{-1}x_2^2.
\label{0v(2.5)}
\end{equation}

Under the assumption of frictionless contact, we have
\begin{equation}
\sigma_{31}({\bf y},0)=0,\quad \sigma_{32}({\bf y},0)=0, \quad {\bf y}=(y_1,y_2)\in\mathbb{R}^2.
\label{0v(2.6)}
\end{equation}
Here the Cartesian coordinate system $(y_1,y_2,z)$ is used such that $y_1=x_1$, $y_2=x_2$, $z=x_3$.

Denoting by $\delta_0$ the indenter's displacement, we formulate the boundary condition on the contact surface as follows:
\begin{equation}
\begin{array}{c}
u_3({\bf y},0)\geq\delta_0-\varphi({\bf y}),\quad
\sigma_{33}({\bf y},0)\leq 0, \\
(u_3({\bf y},0)-\delta_0+\varphi({\bf y}))\sigma_{33}({\bf y},0)=0,\quad 
{\bf y}\in\mathbb{R}^2.
\end{array}
\label{0v(2.7)}
\end{equation}

On the rigid substrate surface (\ref{0v(2.1)}), we will have
\begin{equation}
u_j({\bf y},H({\bf y}))=0,\quad {\bf y}\in\mathbb{R}^2 \quad (j=1,2,3).
\label{0v(2.8)}
\end{equation}

Assuming that the layer is relatively thin in comparison to the characteristic dimensions of $\omega$, we introduce a small dimensionless parameter $\varepsilon$ and set
\begin{equation}
\delta_0=\varepsilon\delta_0^*,\quad R_1=\varepsilon^{-1}R_1^*,\quad R_2=\varepsilon^{-1}R_2^*,
\label{0v(2.9)}
\end{equation}
\begin{equation}
H({\bf y})=\varepsilon h_*(1+\varepsilon\psi_*({\bf y})),
\label{0v(2.10)}
\end{equation}
where $\delta_0^*$, $R_1^*$, and $R_2^*$ are assumed to be comparable with $h_*$, and, moreover,
$\vert\psi_*({\bf y})\vert\leq h_*$ for ${\bf y}\in\mathbb{R}^2$.

The problem is to calculate the contact pressure distribution 
\begin{equation}
p(x_1,x_2)=-\sigma_{33}(x_1,x_2,0),\quad (x_1,x_2)\in\omega.
\label{0v(2.11)}
\end{equation}

An important characteristic of the problem is the contact force required to indent the layer 
\begin{equation}
P=\iint\limits_\omega p({\bf y})\,d{\bf y}.
\label{0v(2.12)}
\end{equation}

Finally, let us introduce the notation
\begin{equation}
h=\varepsilon h_*,\quad \tilde{H}({\bf y})=\varepsilon^2 \tilde{H}_*({\bf y}),
\label{0v(2.13)}
\end{equation}
where 
\begin{equation}
\tilde{H}_*({\bf y})=h_*\psi_*({\bf y}).
\label{0v(2.14)}
\end{equation}

Thus, the following relation takes place:
\begin{equation}
H({\bf y})=h+\tilde{H}({\bf y}).
\label{0v(2.15)}
\end{equation}
Here, $h$ is an average thickness, $\tilde{H}({\bf y})$ is a small variation such that 
$\tilde{H}({\bf y})\ll h$ (see Figure~\ref{figure_strip.pdf}).

\section{Asymptotic Ansatz}
\label{0vSection3}

First, we introduce the so-called stretched coordinate
\begin{equation}
\zeta=\varepsilon^{-1}z.
\label{0v(3.1)}
\end{equation}

Now, substituting (\ref{0v(2.3)}) and (\ref{0v(2.4)}) into Eqs.~(\ref{0v(2.2)}) and taking into account (\ref{0v(3.1)}), we arrive at the following Lam\'e system for the displacement vector 
${\bf u}=({\bf v},w)$:
\begin{equation}
\varepsilon^{-2}\mu\frac{\partial^2{\bf v}}{\partial\zeta^2}+
\varepsilon^{-1}(\lambda+\mu)\nabla_y\frac{\partial w}{\partial\zeta}+
\mu\nabla_y\cdot\nabla_y{\bf v}+(\lambda+\mu)\nabla_y\nabla_y\cdot{\bf v}=0,
\label{0v(3.2)}
\end{equation}
\begin{equation}
\varepsilon^{-2}(2\mu+\lambda)\frac{\partial^2 w}{\partial\zeta^2}+
\varepsilon^{-1}(\lambda+\mu)\nabla_y\cdot\frac{\partial {\bf v}}{\partial\zeta}+
\mu\nabla_y\cdot\nabla_y w=0.
\label{0v(3.3)}
\end{equation}
Here, $\nabla_y=(\partial/\partial y_1,\partial/\partial y_2)$ is the nabla differential operator, and the dot denotes the scalar product, so that $\nabla_y\cdot\nabla_y=\Delta_y$ is the Laplace operator. 

Correspondingly, the boundary condition (\ref{0v(2.6)}) takes the form
\begin{equation}
\varepsilon^{-1}\frac{\partial{\bf v}}{\partial\zeta}+
\nabla_y w\Bigr\vert_{\zeta=0}=0.
\label{0v(3.4)}
\end{equation}

In view of (\ref{0v(2.7)}) and (\ref{0v(2.9)}), we will have
\begin{equation}
w({\bf y},0)=\varepsilon(\delta_0^*-\varphi^*({\bf y})),\quad {\bf y}\in\omega,
\label{0v(3.5)}
\end{equation}
where we introduced the notation (see Eqs.~(\ref{0v(2.5)}) and (\ref{0v(2.9)}))
\begin{equation}
\varphi^*({\bf y})=(2R_1^*)^{-1}y_1^2+(2R_2^*)^{-1}y_2^2.
\label{0v(3.6)}
\end{equation}

Furthermore, by stretching the vertical coordinate, formula (\ref{0v(2.11)}) is transformed as 
\begin{equation}
-p({\bf y})=\varepsilon^{-1}(2\mu+\lambda)\frac{\partial w}{\partial\zeta}+
\lambda\nabla_y\cdot {\bf v}\Bigr\vert_{\zeta=0},\quad {\bf y}\in\omega.
\label{0v(3.7)}
\end{equation}

Finally, the boundary conditions (\ref{0v(2.8)}) on the substrate surface (see Eqs.~(\ref{0v(2.1)}) and (\ref{0v(2.10)}))
\begin{equation}
\zeta=h_*(1+\varepsilon\psi_*({\bf y}))
\label{0v(3.8)}
\end{equation}
take the following form:
\begin{equation}
{\bf v}({\bf y},h_*+\varepsilon h_*\psi_*({\bf y}))=0,\quad
w({\bf y},h_*+\varepsilon h_*\psi_*({\bf y}))=0.
\label{0v(3.9)}
\end{equation}

Observe that among Eqs.~(\ref{0v(3.2)})\,--\,(\ref{0v(3.5)}), (\ref{0v(3.7)}), and (\ref{0v(3.9)}), there are only two inhomogeneous ones, namely, (\ref{0v(3.5)}) and (\ref{0v(3.7)}). So, the form of (\ref{0v(3.5)}) suggests the asymptotic expansion
\begin{equation}
w({\bf y},\zeta)=\varepsilon w^0({\bf y},\zeta)+\varepsilon^2 w^1({\bf y},\zeta)+\ldots\,.
\label{0v(3.10)}
\end{equation}

Now, in view of (\ref{0v(3.7)}) and (\ref{0v(3.10)}), we may suggest that 
$p({\bf y})=O(1)$ as $\varepsilon\to 0$. However, taking into account the homogeneous conditions (\ref{0v(3.4)}) and (\ref{0v(3.9)}), we put
\begin{equation}
{\bf v}({\bf y},\zeta)=\varepsilon {\bf v}^0({\bf y},\zeta)+\varepsilon^2 {\bf v}^1({\bf y},\zeta)+\ldots\,.
\label{0v(3.11)}
\end{equation}

We emphasize that the asymptotic Ansatz (\ref{0v(3.10)}), (\ref{0v(3.11)}) is valid only inside the contact region~$\omega$. In other words, a plane boundary layer should be constructed near the edge of the contact area. We refer to \cite{Argatov2005,Chadwick2002} for more details.

\section{Derivation of asymptotic expansions}
\label{0vSection4}

Substitution of (\ref{0v(3.10)}) and (\ref{0v(3.11)}) into Eqs.~(\ref{0v(3.2)}) and (\ref{0v(3.3)}) gives
\begin{equation}
\begin{array}{l}
\displaystyle
\varepsilon^{-2}\mu\frac{\partial^2{\bf v}^0}{\partial\zeta^2}+
\varepsilon^{-1}\biggl((\lambda+\mu)\nabla_y\frac{\partial w^0}{\partial\zeta}+
\mu\frac{\partial^2{\bf v}^1}{\partial\zeta^2}\biggr) \\
\displaystyle
{}+\varepsilon^0\biggl(
\mu\Delta_y{\bf v}^0+(\lambda+\mu)\nabla_y\nabla_y\cdot{\bf v}^0
+(\lambda+\mu)\nabla_y\frac{\partial w^1}{\partial\zeta}+
\mu\frac{\partial^2{\bf v}^2}{\partial\zeta^2}\biggr)+\ldots=0,
\end{array}
\label{0v(4.1)}
\end{equation}
\begin{equation}
\begin{array}{l}
\displaystyle
\varepsilon^{-2}(2\mu+\lambda)\frac{\partial^2 w^0}{\partial\zeta^2}+
\varepsilon^{-1}\biggl((2\mu+\lambda)\frac{\partial^2 w^1}{\partial\zeta^2}+
(\lambda+\mu)\nabla_y\cdot\frac{\partial {\bf v}^0}{\partial\zeta}\biggr) \\
\displaystyle
{}+\varepsilon^0\biggl((2\mu+\lambda)\frac{\partial^2 w^2}{\partial\zeta^2}+
(\lambda+\mu)\nabla_y\cdot\frac{\partial {\bf v}^1}{\partial\zeta}
+\mu\Delta_y w^0\biggr)+\ldots=0.
\end{array}
\label{0v(4.2)}
\end{equation}

Further, the substitution of (\ref{0v(3.10)}) and (\ref{0v(3.11)}) into the boundary conditions (\ref{0v(3.4)}) and (\ref{0v(3.7)}) at the contact region yields
\begin{equation}
\varepsilon^{-1}\frac{\partial{\bf v}^0}{\partial\zeta}+
\varepsilon^0\Bigl(\nabla_y w^0+\frac{\partial{\bf v}^1}{\partial\zeta}\Bigr)
+\varepsilon\Bigl(\nabla_y w^1+\frac{\partial{\bf v}^2}{\partial\zeta}\Bigr)+\ldots
\Bigr\vert_{\zeta=0}=0,
\label{0v(4.3)}
\end{equation}
\begin{equation}
\begin{array}{l}
\displaystyle
(2\mu+\lambda)\frac{\partial w^0}{\partial\zeta}+
\varepsilon\Bigl(\lambda\nabla_y\cdot {\bf v}^0
+(2\mu+\lambda)\frac{\partial w^1}{\partial\zeta}\Bigr) \\
\displaystyle
{}+\varepsilon^2\Bigl(\lambda\nabla_y\cdot {\bf v}^1
+(2\mu+\lambda)\frac{\partial w^2}{\partial\zeta}\Bigr)+\ldots
\Bigr\vert_{\zeta=0}=-p({\bf y}).
\end{array}
\label{0v(4.4)}
\end{equation}

Finally, the substitution of (\ref{0v(3.10)}) and (\ref{0v(3.11)}) into (\ref{0v(3.9)}) leads to the boundary conditions
\begin{equation}
{\bf v}^0+\varepsilon\Bigl({\bf v}^1+\tilde{H}_*\frac{\partial{\bf v}^0}{\partial\zeta}\Bigr)
+\varepsilon^2\Bigl({\bf v}^2+\tilde{H}_*\frac{\partial{\bf v}^1}{\partial\zeta}
+\frac{\tilde{H}_*^2}{2}\frac{\partial^2{\bf v}^0}{\partial\zeta^2}\Bigr)
+\ldots\Bigr\vert_{\zeta=h_*}=0,
\label{0v(4.5)}
\end{equation}
\begin{equation}
w^0+\varepsilon\Bigl(w^1+\tilde{H}_*\frac{\partial w^0}{\partial\zeta}\Bigr)
+\varepsilon^2\Bigl(w^2+\tilde{H}_*\frac{\partial w^1}{\partial\zeta}
+\frac{\tilde{H}_*^2}{2}\frac{\partial^2 w^0}{\partial\zeta^2}\Bigr)
+\ldots\Bigr\vert_{\zeta=h_*}=0,
\label{0v(4.6)}
\end{equation}
where the notation (\ref{0v(2.14)}) was taken into account.

Thus, on the basis of Eqs.~(\ref{0v(4.1)})\,--\,(\ref{0v(4.6)}), we arrive at a recurrence system of boundary-value problems for the functions ${\bf v}^k$ and $w^k$ ($k=0,1,\ldots$). In the next two sections, we will construct the first several terms of the asymptotic series (\ref{0v(3.10)}) and (\ref{0v(3.11)}).

\section{Asymptotic model for a compressible elastic layer}
\label{0vSection5}

According to (\ref{0v(4.1)})\,--\,(\ref{0v(4.6)}), the first-order problem takes the form 
\begin{equation}
(2\mu+\lambda)\frac{\partial^2 w^0}{\partial\zeta^2}=0, \quad \zeta\in(0,h_*),\quad
(2\mu+\lambda)\frac{\partial w^0}{\partial\zeta}\Bigr\vert_{\zeta=0}=-p({\bf y}),\quad
w^0\bigr\vert_{\zeta=h_*}=0;
\label{0v(5.1)}
\end{equation}
\begin{equation}
\mu\frac{\partial^2{\bf v}^0}{\partial\zeta^2}=0, \quad \zeta\in(0,h_*),\quad
\frac{\partial{\bf v}^0}{\partial\zeta}\Bigr\vert_{\zeta=0}=0,\quad
{\bf v}^0\bigr\vert_{\zeta=h_*}=0.
\label{0v(5.2)}
\end{equation}

From (\ref{0v(5.1)}) and (\ref{0v(5.2)}), it immediately follows that 
\begin{equation}
w^0({\bf y},\zeta)=\frac{p({\bf y})}{2\mu+\lambda}(h_*-\zeta),
\label{0v(5.3)}
\end{equation}
\begin{equation}
{\bf v}^0({\bf y},\zeta)\equiv 0.
\label{0v(5.4)}
\end{equation}

In view of (\ref{0v(5.4)}), the second-order problem, derived from Eqs.~(\ref{0v(4.1)})\,--\,(\ref{0v(4.6)}), takes the following form:
\begin{equation}
\frac{\partial^2 w^1}{\partial\zeta^2}=0, \quad \zeta\in(0,h_*),\quad
\frac{\partial w^1}{\partial\zeta}\Bigr\vert_{\zeta=0}=0,\quad
w^1\bigr\vert_{\zeta=h_*}=-\tilde{H}_*({\bf y})\frac{\partial w^0}{\partial\zeta}\Bigr\vert_{\zeta=h_*};
\label{0v(5.5)}
\end{equation}
\begin{equation}
\mu\frac{\partial^2{\bf v}^1}{\partial\zeta^2}=-(\lambda+\mu)\nabla_y\frac{\partial w^0}{\partial\zeta}, \quad \zeta\in(0,h_*),\quad
\frac{\partial{\bf v}^1}{\partial\zeta}\Bigr\vert_{\zeta=0}=-\nabla_y w^0\bigr\vert_{\zeta=0},\quad
{\bf v}^1\bigr\vert_{\zeta=h_*}=0.
\label{0v(5.6)}
\end{equation}

It can be easily shown that the solution of the problem (\ref{0v(5.5)}) reads as
\begin{equation}
w^1({\bf y},\zeta)=\frac{p({\bf y})}{2\mu+\lambda}\tilde{H}_*({\bf y}).
\label{0v(5.7)}
\end{equation}

On the other hand, the unique solution of the problem (\ref{0v(5.6)}) can be represented as
\begin{equation}
{\bf v}^1({\bf y},\zeta)=\Psi(\zeta)\nabla_y p({\bf y}),
\label{0v(5.8)}
\end{equation}
where we introduced the notation 
\begin{equation}
\Psi(\zeta)=-\frac{\lambda+\mu}{2\mu(2\mu+\lambda)}(h_*^2-\zeta^2)
+\frac{h_*}{2\mu+\lambda}(h_*-\zeta).
\label{0v(5.9)}
\end{equation}

Thus, collecting Eqs.~(\ref{0v(3.10)}), (\ref{0v(5.3)}), and (\ref{0v(5.7)}), we obtain the two-term asymptotic approximation for the normal displacement
\begin{equation}
w({\bf y},\zeta)\simeq \varepsilon\frac{p({\bf y})}{2\mu+\lambda}(h_*-\zeta)
+\varepsilon^2\frac{p({\bf y})}{2\mu+\lambda}\tilde{H}_*({\bf y}).
\label{0v(5.10)}
\end{equation}
By taking into account the scaling relations (\ref{0v(2.13)}), we rewrite (\ref{0v(5.10)}) in the form
\begin{equation}
u_3({\bf y},z)\simeq \frac{p({\bf y})}{2\mu+\lambda}(h-z)
+\frac{p({\bf y})}{2\mu+\lambda}\tilde{H}({\bf y}).
\label{0v(5.11)}
\end{equation}

Now, substituting the expression (\ref{0v(5.11)}) into the contact condition 
\begin{equation}
u_3({\bf y},0)=\delta_0-\varphi({\bf y}),\quad{\bf y}\in\omega,
\label{0v(5.12)}
\end{equation}
we derive the following equation for the contact pressure density:
\begin{equation}
\frac{h+\tilde{H}({\bf y})}{2\mu+\lambda}p({\bf y})=\delta_0-\varphi({\bf y}),\quad{\bf y}\in\omega.
\label{0v(5.13)}
\end{equation}

In view of the condition $p({\bf y})>0$ for ${\bf y}\in\omega$, we get
\begin{equation}
p({\bf y})=\frac{2\mu+\lambda}{h+\tilde{H}({\bf y})}
\bigl(\delta_0-\varphi({\bf y})\bigr)_+,
\label{0v(5.14)}
\end{equation}
where $(x)_+=\max\{x,0\}$ is the positive-part function.

Now, invoking the notation (\ref{0v(2.15)}) for the variable thickness of the elastic layer, we rewrite (\ref{0v(5.14)}) as follows:
\begin{equation}
p({\bf y})=\frac{2\mu+\lambda}{H({\bf y})}\bigl(\delta_0-\varphi({\bf y})\bigr)_+.
\label{0v(5.15)}
\end{equation}

Formula (\ref{0v(5.15)}) shows that a thin compressible elastic layer deforms like a Winkler foundation with the variable foundation modulus 
\begin{equation}
k({\bf y})=\frac{2\mu+\lambda}{H({\bf y})}.
\label{0v(5.16)}
\end{equation}

Finally, let us recall that the Lam\'e parameters $\lambda$ and $\mu$ are related to Young's modulus, $E$, and Poisson's ratio, $\nu$, by formulas
\begin{equation}
\lambda=\frac{E\nu}{(1+\nu)(1-2\nu)},\quad \mu=\frac{E}{2(1+\nu)}.
\label{0v(5.17)}
\end{equation}

Now, in view of (\ref{0v(5.17)}), it is readily seen from (\ref{0v(5.16)}) that $k({\bf y})\to\infty$ as $\nu\to 0{.}5$. This implies that the case of an incompressible elastic layer with $\nu=0{.}5$ requires a special consideration. 

\section{Asymptotic model for an incompressible elastic layer}
\label{0vSection6}

Let us continue the process of constructing terms of the asymptotic expansions (\ref{0v(3.10)}) and (\ref{0v(3.11)}). In view of (\ref{0v(5.7)}), Eqs.~(\ref{0v(4.1)})\,--\,(\ref{0v(4.6)}) yield the third-order problem
$$
(2\mu+\lambda)\frac{\partial^2 w^2}{\partial\zeta^2}=
-(\lambda+\mu)\nabla_y\cdot\frac{\partial {\bf v}^1}{\partial\zeta}
-\mu\Delta_y w^0, \quad \zeta\in(0,h_*),
$$
\begin{equation}
(2\mu+\lambda)\frac{\partial w^2}{\partial\zeta}\Bigr\vert_{\zeta=0}=
-\lambda\nabla_y\cdot {\bf v}^1\bigr\vert_{\zeta=0},\quad
w^2\bigr\vert_{\zeta=h_*}=0;
\label{0v(6.1)}
\end{equation}
\begin{equation}
\frac{\partial^2{\bf v}^2}{\partial\zeta^2}=0, \quad \zeta\in(0,h_*),\quad
\frac{\partial{\bf v}^2}{\partial\zeta}\Bigr\vert_{\zeta=0}=
-\nabla_y w^1\bigr\vert_{\zeta=0},\quad
{\bf v}^2\bigr\vert_{\zeta=h_*}=
-\tilde{H}_*({\bf y})\frac{\partial {\bf v}^1}{\partial\zeta}\Bigr\vert_{\zeta=h_*}.
\label{0v(6.2)}
\end{equation}

Substituting (\ref{0v(5.3)}) and (\ref{0v(5.8)}) into Eqs.~(\ref{0v(6.1)}), we derive the boundary-value problem
$$
\frac{\partial^2 w^2}{\partial\zeta^2}=-\frac{\lambda\Delta_y p({\bf y})}{\mu(2\mu+\lambda)^2}
[(2\mu+\lambda)\zeta-\mu h_*], \quad \zeta\in(0,h_*),
$$
\begin{equation}
\frac{\partial w^2}{\partial\zeta}\Bigr\vert_{\zeta=0}=
-\frac{\lambda(\mu-\lambda)h_*^2}{2\mu(2\mu+\lambda)^2}\Delta_y p({\bf y}),\quad
w^2\bigr\vert_{\zeta=h_*}=0.
\label{0v(6.3)}
\end{equation}

It can be checked that the solution to (\ref{0v(6.3)}) can be condensed to the form
\begin{eqnarray}
w^2({\bf y},\zeta) & = & \frac{\Delta_y p({\bf y})}{6\mu(2\mu+\lambda)^2}\bigl\{
3\lambda\mu h_*(\zeta^2-h_*^2)\nonumber\\
{} & {} & {}-\lambda(2\mu+\lambda)(\zeta^3-h_*^3)
-3\lambda(\mu-\lambda)h_*^2(\zeta-h_*)\bigr\}.
\label{0v(6.4)}
\end{eqnarray}

Further, in view of (\ref{0v(5.7)}) and (\ref{0v(5.8)}), the problem (\ref{0v(6.2)}) takes the form
\begin{equation}
\frac{\partial^2{\bf v}^2}{\partial\zeta^2}=0, \quad \zeta\in(0,h_*),\quad
\frac{\partial{\bf v}^2}{\partial\zeta}\Bigr\vert_{\zeta=0}=
-\frac{\nabla_y (p\tilde{H}_*)}{2\mu+\lambda},\quad
{\bf v}^2\bigr\vert_{\zeta=h_*}=
-\frac{\lambda h_*}{\mu(2\mu+\lambda)}\tilde{H}_*\nabla_y p.
\label{0v(6.5)}
\end{equation}
Here the arguments of functions $p({\bf y})$ and $\tilde{H}_*({\bf y})$ are omitted for clarity. 

It can be easily verified that the solution to (\ref{0v(6.5)}) has the form
\begin{equation}
{\bf v}^2({\bf y},\zeta)=
\frac{h_*-\zeta}{2\mu+\lambda}\nabla_y (p\tilde{H}_*)
-\frac{\lambda h_*}{\mu(2\mu+\lambda)}\tilde{H}_*\nabla_y p.
\label{0v(6.6)}
\end{equation}

We emphasize that in contrast to the first two term approximation (\ref{0v(5.10)}), the third term (\ref{0v(6.4)}) does not vanish at the contact surface in the limit as $\nu\to 0{.}5$. Indeed, formula (\ref{0v(6.4)}) yields
\begin{equation}
w^2({\bf y},0)=-\frac{h_*^3 \lambda(\lambda-\mu)}{3\mu(2\mu+\lambda)^2}\Delta_y p({\bf y}),
\label{0v(6.7)}
\end{equation}
where, in view of (\ref{0v(5.17)}), as $\nu\to 0{.}5$, we will have
\begin{equation}
\frac{\lambda(\lambda-\mu)}{\mu(2\mu+\lambda)^2}=
\frac{\nu(1+\nu)(4\nu-1)}{E(1-\nu)^2}\rightarrow\frac{3}{E}.
\label{0v(6.7a)}
\end{equation}

Finally, in order to construct a correction for the leading asymptotic term (\ref{0v(6.4)}), we consider the following problem:
$$
(2\mu+\lambda)\frac{\partial^2 w^3}{\partial\zeta^2}=
-(\lambda+\mu)\nabla_y\cdot\frac{\partial {\bf v}^2}{\partial\zeta}
-\mu\Delta_y w^1, \quad \zeta\in(0,h_*),
$$
\begin{equation}
(2\mu+\lambda)\frac{\partial w^3}{\partial\zeta}\Bigr\vert_{\zeta=0}=
-\lambda\nabla_y\cdot {\bf v}^2\bigr\vert_{\zeta=0},\quad
w^3\bigr\vert_{\zeta=h_*}=
-\tilde{H}_*({\bf y})\frac{\partial w^2}{\partial\zeta}\Bigr\vert_{\zeta=h_*}.
\label{0v(6.8)}
\end{equation}

Substituting the expressions (\ref{0v(5.7)}), (\ref{0v(6.4)}), and (\ref{0v(6.6)}) into Eqs.~(\ref{0v(6.8)}), we get
\begin{equation}
\frac{\partial^2 w^3}{\partial\zeta^2}=
\frac{\lambda}{(2\mu+\lambda)^2}\Delta_y (p\tilde{H}_*), \quad \zeta\in(0,h_*),
\label{0v(6.9a)}
\end{equation}
\begin{equation}
\frac{\partial w^3}{\partial\zeta}\Bigr\vert_{\zeta=0}=
\frac{\lambda h_*}{\mu(2\mu+\lambda)}\bigl(
\lambda\nabla_y\cdot(\tilde{H}_*\nabla_y p)-
\mu \Delta_y (p\tilde{H}_*)\bigr),\quad
w^3\bigr\vert_{\zeta=h_*}=
\frac{\lambda h_*^2}{2(2\mu+\lambda)^2}\tilde{H}_*\Delta_y p.
\label{0v(6.9)}
\end{equation}

Integrating Eq.~(\ref{0v(6.9a)}), we obtain 
\begin{equation}
w^3({\bf y},0)=\frac{\lambda}{2(2\mu+\lambda)^2}\Delta_y (p\tilde{H}_*)\zeta^2+
C_1({\bf y})\zeta+C_0({\bf y}),
\label{0v(6.10)}
\end{equation}
where the integration functions $C_1({\bf y})$ and $C_0({\bf y})$ are restricted by the boundary conditions 
(\ref{0v(6.9)}). So, it can be checked that
\begin{equation}
C_1({\bf y})=\frac{\lambda h_*}{\mu(2\mu+\lambda)}\bigl(
\lambda\nabla_y\cdot(\tilde{H}_*\nabla_y p)-
\mu \Delta_y (p\tilde{H}_*)\bigr),
\label{0v(6.11)}
\end{equation}
\begin{equation}
C_0({\bf y})=\frac{\lambda h_*^2}{2(2\mu+\lambda)^2}
[\tilde{H}_*\Delta_y p+\Delta_y (p\tilde{H}_*)]
-\frac{\lambda^2 h_*^2}{\mu(2\mu+\lambda)^2}\nabla_y\cdot(\tilde{H}_*\nabla_y p).
\label{0v(6.12)}
\end{equation}

From (\ref{0v(6.10)}), it immediately follows that 
\begin{equation}
w^3({\bf y},0)=C_0({\bf y}).
\label{0v(6.13)}
\end{equation}

Note that the elastic constants entering Eq.~(\ref{0v(6.12)}) are evaluated in terms of the engineering elastic constants as 
\begin{equation}
\frac{\lambda}{(2\mu+\lambda)^2}=\frac{\nu(1+\nu)(1-2\nu)}{E(1-\nu)^2},\quad
\frac{\lambda^2}{\mu(2\mu+\lambda)^2}=\frac{2\nu^2(1+\nu)}{E(1-\nu)^2}.
\label{0v(6.14)}
\end{equation}
Hence, in view of (\ref{0v(6.14)}), the first term in (\ref{0v(6.12)}) disappears as $\nu\to 0{.}5$, and we get
\begin{equation}
C_0({\bf y})\bigr\vert_{\nu=0{.}5}=-\frac{3h_*^2}{E}\nabla_y\cdot(\tilde{H}_*\nabla_y p).
\label{0v(6.15)}
\end{equation}

Thus, collecting Eqs.~(\ref{0v(3.10)}), (\ref{0v(6.7)}), (\ref{0v(6.13)}), and (\ref{0v(6.15)}), we obtain the following two-term asymptotic approximation for the normal displacement at the contact surface in the case of incompressible elastic layer:
\begin{equation}
w({\bf y},0)\simeq -\varepsilon^3\frac{h_*^3}{E}\Delta_y p({\bf y})
-\varepsilon^4\frac{3h_*^2}{E}\nabla_y\cdot(\tilde{H}_*({\bf y})\nabla_y p({\bf y})).
\label{0v(6.16)}
\end{equation}
Recollecting the scaling relations (\ref{0v(2.13)}), we rewrite (\ref{0v(6.16)}) in the form
\begin{equation}
u_3({\bf y},0)\simeq -\frac{h^3}{E}\Delta_y p({\bf y})
-\frac{3h^2}{E}\nabla_y\cdot(\tilde{H}({\bf y})\nabla_y p({\bf y})).
\label{0v(6.17)}
\end{equation}

Now, substituting the expression (\ref{0v(6.17)}) into the contact condition (\ref{0v(5.12)}), we arrive at a partial differential equation in the domain $\omega$ with respect to the function $p({\bf y})$. According to the asymptotic analysis \cite{Chadwick2002}, at the contour $\Gamma$ of $\omega$, we impose the following boundary conditions:
\begin{equation}
p({\bf y})=0,\quad \frac{\partial p}{\partial n}({\bf y})=0,\quad {\bf y}\in\Gamma.
\label{0v(6.18)}
\end{equation}
Here, $\partial/\partial n$ is the normal derivative. We stress that the location of the contour $\Gamma$ must be determined as part of the solution. 

\section{Comparison of the obtained result with the solution for the 2D case}
\label{0vSection7}

Collecting Eqs.~(\ref{0v(3.10)}), (\ref{0v(5.3)}), (\ref{0v(5.7)}), (\ref{0v(6.4)}), and (\ref{0v(6.10)}), we obtain
\begin{eqnarray}
w({\bf y},0) & \simeq & \varepsilon\frac{h_*}{2\mu+\lambda}p({\bf y})
+\varepsilon^2\frac{\tilde{H}_*({\bf y})}{2\mu+\lambda}p({\bf y})
-\varepsilon^3\frac{h_*^3\lambda(\lambda-\mu)}{3\mu(2\mu+\lambda)^2}\Delta_y p({\bf y})
\nonumber\\
{} & {} & {}+\varepsilon^4\biggl\{\frac{\lambda h_*^2}{2(2\mu+\lambda)^2}
\bigl[\tilde{H}_*({\bf y})\Delta_y p({\bf y})+\Delta_y \bigl(p({\bf y})\tilde{H}_*({\bf y})\bigr)\bigr]
\nonumber\\
{} & {} & {}-\frac{\lambda^2 h_*^2}{\mu(2\mu+\lambda)^2}\nabla_y\cdot\bigl(\tilde{H}_*({\bf y})\nabla_y p({\bf y})\bigr)\biggr\}.
\label{0v(7.1)}
\end{eqnarray}

Substituting the asymptotic expansion (\ref{0v(7.1)}) into the contact condition (\ref{0v(3.5)}) and using the notation (\ref{0v(2.14)}), we derive the following equation for the contact pressure density:
\begin{eqnarray}
{} & {} & p({\bf y})+\varepsilon \psi_*({\bf y})p({\bf y})
-\varepsilon^2\frac{h_*^3\lambda(\lambda-\mu)}{3\mu(2\mu+\lambda)}\Delta_y p({\bf y})
\nonumber\\
{} & {} & {}+\varepsilon^3\frac{\lambda h_*^2}{2\mu(2\mu+\lambda)}
\biggl\{\mu \bigl[\psi_*({\bf y})\Delta_y p({\bf y})+\Delta_y \bigl(p({\bf y})\psi_*({\bf y})\bigr)\bigr]
\nonumber\\
{} & {} & {}-2\lambda\nabla_y\cdot\bigl(\psi_*({\bf y})\nabla_y p({\bf y})\bigr)\biggr\}
=\frac{2\mu+\lambda}{h_*}f^*({\bf y}).
\label{0v(7.2)}
\end{eqnarray}
Here we also introduced a shorthand notation for the right-hand side of (\ref{0v(3.5)}), i.\,e.,
\begin{equation}
f^*({\bf y})=\delta_0^*-\varphi^*({\bf y}).
\label{0v(7.3)}
\end{equation}

It should be noted that Eq.~(\ref{0v(7.2)}) is applied for the case of compressible materials when its right-hand side makes sense.

By applying a perturbation method, a solution to Eq.~(\ref{0v(7.2)}) is represented in the form
\begin{equation}
p({\bf y})\simeq \frac{2\mu+\lambda}{h_*}\bigl(
\sigma_0({\bf y})+\varepsilon\sigma_1({\bf y})+\varepsilon^2\sigma_2({\bf y})+\varepsilon^3\sigma_3({\bf y})
\bigr).
\label{0v(7.4)}
\end{equation}

After the substitution of (\ref{0v(7.4)}) into (\ref{0v(7.2)}), we straightforwardly obtain 
\begin{equation}
\sigma_0=f^*,\quad \sigma_1=-\psi_* f^*,\quad \sigma_2=\psi_*^2 f^*
+\frac{h_*^2\lambda(\lambda-\mu)}{3\mu(2\mu+\lambda)}\Delta_y f^*,
\label{0v(7.5a)}
\end{equation}
\begin{eqnarray}
\sigma_3 & = & -\psi_*^3 f^*
-\frac{h_*^2\lambda(\lambda-\mu)}{3\mu(2\mu+\lambda)}
(\psi_*\Delta_y f^*+\Delta_y(\psi_* f^*))
\nonumber\\
{} & {} & {}-\frac{h_*^2\lambda}{2\mu(2\mu+\lambda)}
\bigl(
\mu [\psi_*\Delta_y f^*+\Delta_y (f^*\psi_*)]
-2\lambda\nabla_y\cdot(\psi_*\nabla_y f^*)\bigr),
\label{0v(7.5)}
\end{eqnarray}
where for the sake of brevity, the argument ${\bf y}$ is omitted. 

Further, making use of the differential identities
$$
\nabla\cdot(\psi\nabla f)=\nabla\psi\cdot\nabla f+\psi\Delta f,
$$
$$
\Delta (f\psi)=\psi\Delta f+f\Delta \psi+2\nabla f\cdot\nabla \psi,
$$
we simplify formula (\ref{0v(7.5)}) as follows:
\begin{equation}
\sigma_3 = -\psi_*^3 f^*
+\frac{h_*^2\lambda(\lambda-\mu)}{3\mu(2\mu+\lambda)}
(\nabla_y f^*\cdot\nabla_y \psi_*+\psi_*\Delta_y f^*)
-\frac{h_*^2\lambda(2\lambda+\mu)}{6\mu(2\mu+\lambda)} f^*\Delta_y \psi_*.
\label{0v(7.6)}
\end{equation}

Now it can be easily checked that the four-term asymptotic expansion (\ref{0v(7.4)}) with the coefficients given by (\ref{0v(7.5a)}) and (\ref{0v(7.6)}) in the 2D case recovers the corresponding solution obtained in 
\cite{VorovichPeninin1971}, where also the next asymptotic term in (\ref{0v(7.4)}) was explicitly written out. 

Finally, observe that formula (\ref{0v(6.12)}) can be transformed into the following one:
\begin{equation}
C_0({\bf y})=-\frac{h_*^2\lambda(\lambda-\mu)}{\mu(2\mu+\lambda)^2}
[\nabla_y \tilde{H}_*\cdot \nabla_y p+\tilde{H}_*\Delta_y p]
+\frac{h_*^2\lambda}{2(2\mu+\lambda)^2}p\Delta_y\tilde{H}_*.
\label{0v(6.12a)}
\end{equation}
Note also that the expression in the brackets in (\ref{0v(6.12a)}) is equal to 
$\nabla_y\cdot(\tilde{H}_* \nabla_y p)$.

\section{Application to sensitivity analysis of articular contact mechanics}
\label{0vSection8}

According to (\ref{0v(6.17)}), (\ref{0v(6.18)}) the refined asymptotic model for contact interaction of thin incompressible layers bonded to rigid substrates looks as follows:
\begin{equation}
-m^{-1}\Delta_y p({\bf y})
-\sum_{\alpha=1}^2\frac{3h_\alpha^2}{E_\alpha}\nabla_y\cdot(\tilde{H}_\alpha({\bf y})\nabla_y p({\bf y}))
=\delta_0-\varphi({\bf y}),\quad {\bf y}\in\tilde\omega,
\label{0v(8.1)}
\end{equation}
\begin{equation}
p({\bf y})=0,\quad \frac{\partial p}{\partial n}({\bf y})=0,\quad {\bf y}\in\tilde\Gamma.
\label{0v(8.2)}
\end{equation}
Here, $\tilde\Gamma$ is the contour of the contact region $\tilde\omega$, and the notation (\ref{0v(1.6)}) was used. 

Let us put
\begin{equation}
p({\bf y})=\bar{p}({\bf y})+\tilde{p}({\bf y}),
\label{0v(8.3)}
\end{equation}
where $\bar{p}({\bf y})$ is the solution to the original asymptotic model (\ref{0v(1.1)}), (\ref{0v(1.2)}).

Then, under the assumption that the thickness variation functions $\tilde{H}_1({\bf y})$ and $\tilde{H}_2({\bf y})$ introduce a small variation into the elliptical contact region $\omega$ corresponding to the density $\bar{p}({\bf y})$, we derive from (\ref{0v(8.1)})\,--\,(\ref{0v(8.3)}) the following limit problem for the variation of the contact pressure density:
\begin{equation}
-m^{-1}\Delta_y \tilde{p}({\bf y})=
\sum_{\alpha=1}^2\frac{3h_\alpha^2}{E_\alpha}\nabla_y\cdot(\tilde{H}_\alpha({\bf y})\nabla_y 
\bar{p}({\bf y})),\quad {\bf y}\in\omega,
\label{0v(8.4)}
\end{equation}
\begin{equation}
\tilde{p}({\bf y})=0,\quad {\bf y}\in\Gamma.
\label{0v(8.5)}
\end{equation}
Here, $\Gamma$ is the contour corresponding to the contact pressure (\ref{0v(1.4)}).

Moreover, the thickness variations $\tilde{H}_1({\bf y})$ and $\tilde{H}_2({\bf y})$  will not greatly influence the resulting force-displacement relationship, if 
\begin{equation}
\iint\limits_\omega \tilde{p}({\bf y})\,d{\bf y}=0.
\label{0v(8.6)}
\end{equation}

Let us derive the conditions for $\tilde{H}_1({\bf y})$ and $\tilde{H}_2({\bf y})$ under which the equality (\ref{0v(8.6)}) holds true. With this aim we consider an auxiliary problem 
\begin{equation}
\Delta_y \Theta({\bf y})=1,\quad {\bf y}\in\omega,\quad
\Theta({\bf y})=0,\quad {\bf y}\in\Gamma
\label{0v(8.Th1)}
\end{equation}
with the solution 
$\Theta({\bf y})=-a_1^2 a_2^2(2(a_1^2+a_2^2))^{-1}\theta({\bf y})$, where
\begin{equation}
\theta({\bf y})=1-\frac{y_1^2}{a_1^2}-\frac{y_2^2}{a_2^2}.
\label{0v(8.Th2)}
\end{equation}

In view of (\ref{0v(8.Th1)}), we rewrite Eq.~(\ref{0v(8.6)}) as
\begin{equation}
\iint\limits_\omega \tilde{p}({\bf y})\Delta_y \Theta({\bf y})\,d{\bf y}=0.
\label{0v(8.7)}
\end{equation}

Now, applying the second Green's formula and taking into account Eqs.~(\ref{0v(8.4)}), (\ref{0v(8.5)}), and (\ref{0v(8.Th1)}), we reduce Eq.~(\ref{0v(8.7)}) to the following one:
\begin{equation}
\iint\limits_\omega \theta({\bf y})\sum_{\alpha=1}^2\frac{h_\alpha^2}{E_\alpha}\nabla_y\cdot(\tilde{H}_\alpha({\bf y})\nabla_y \bar{p}({\bf y}))\,d{\bf y}=0.
\label{0v(8.8)}
\end{equation}

After rewriting Eq.~(\ref{0v(8.8)}) in the form
$$
\iint\limits_\omega \theta({\bf y})\sum_{\beta=1}^2\frac{\partial}{\partial y_\beta}
\biggl(\frac{\partial\bar{p}}{\partial y_\beta}({\bf y})
\sum_{\alpha=1}^2\frac{h_\alpha^2}{E_\alpha}\tilde{H}_\alpha({\bf y})\biggr)d{\bf y}=0
$$
and integrating by parts with (\ref{0v(8.Th2)}) taken into account, we get
\begin{eqnarray}
{} & {} & -\iint\limits_\omega\sum_{\beta=1}^2 
\frac{2y_\beta}{a_\beta^2}
\frac{\partial \bar{p}}{\partial y_\beta}({\bf y})
\sum_{\alpha=1}^2\frac{h_\alpha^2}{E_\alpha}\tilde{H}_\alpha({\bf y})\,d{\bf y}
\nonumber\\
{} & {} & {}
+\int\limits_\Gamma \theta({\bf y})\sum_{\beta=1}^2\cos(n,y_\beta)
\frac{\partial \bar{p}}{\partial y_\beta}({\bf y})
\sum_{\alpha=1}^2\frac{h_\alpha^2}{E_\alpha}\tilde{H}_\alpha({\bf y})\,ds_y=0.
\label{0v(8.9)}
\end{eqnarray}

It is clear that the line integral in (\ref{0v(8.9)}) vanishes due to the boundary condition (\ref{0v(8.Th2)}). Hence, taking into account the exact expression (\ref{0v(1.4)}) for $\bar{p}({\bf y})$, we finally transform Eq.~(\ref{0v(8.9)}) into the following one:
\begin{equation}
\sum_{\alpha=1}^2\frac{h_\alpha^2}{E_\alpha}\iint\limits_\omega
\tilde{H}_\alpha({\bf y})\rho({\bf y})\,d{\bf y}=0.
\label{0v(8.10)}
\end{equation}
Here we introduced the notation
\begin{equation}
\rho({\bf y})=\biggl(\frac{s y_1^2}{a_1^2}+\frac{y_2^2}{s a_2^2}\biggr)
\biggl(1-\frac{y_1^2}{a_1^2}-\frac{y_2^2}{a_2^2}\biggr).
\label{0v(8.11)}
\end{equation}

Based on the derived Eq.~(\ref{0v(8.10)}), we suggest the following optimization criterion for determining the average thicknesses $h_1$ and $h_2$:
\begin{equation}
\min_{h_\alpha}\iint\limits_{\omega_*}
(H_\alpha({\bf y})-h_\alpha)^2\rho_*({\bf y})\,d{\bf y}.
\label{0v(8.12)}
\end{equation}
Here, $\omega_*$ is an characteristic elliptic domain with semi-axes $a_*$ and $b_*$. In particular, in the capacity of $\omega_*$ one can take the average contact area for a class of admissible contact loadings, and $\rho_*({\bf y})$ is given by
\begin{equation}
\rho_*({\bf y})=\biggl(\frac{s^* y_1^2}{{a_1^*}^2}+\frac{y_2^2}{s^* {a_2^*}^2}\biggr)
\biggl(1-\frac{y_1^2}{{a_1^*}^2}-\frac{y_2^2}{{a_2^*}^2}\biggr),
\label{0v(8.11b)}
\end{equation}
where $s^*=a_2^*/a_1^*$ is the aspect ratio of $\omega_*$.

It is clear that the necessary optimality condition for (\ref{0v(8.12)}) has the form
\begin{equation}
\iint\limits_{\omega_*}
(H_\alpha({\bf y})-h_\alpha)\rho_*({\bf y})\,d{\bf y}=0,
\label{0v(8.13)}
\end{equation}
from where it follows that 
\begin{equation}
h_\alpha=\biggl(\iint\limits_{\omega_*}\rho_*({\bf y})\,d{\bf y}\biggr)^{-1}
\iint\limits_{\omega_*}
H_\alpha({\bf y})\rho_*({\bf y})\,d{\bf y}.
\label{0v(8.14)}
\end{equation}

It is left to show that Eq.~(\ref{0v(8.10)}) follows from (\ref{0v(8.13)}) if $\omega_*$ coincides with $\omega$. Indeed, in view of (\ref{0v(2.15)}), Eq.~(\ref{0v(8.13)}) is equivalent to the following one:
\begin{equation}
\iint\limits_{\omega_*}
\tilde{H}_\alpha({\bf y})\rho_*({\bf y})\,d{\bf y}=0,\quad \alpha=1,2.
\label{0v(8.15H)}
\end{equation}
Now, adding the two equations above multiplied by $E_\alpha^{-1}h_\alpha^2$, $\alpha=1,2$, respectively, we arrive at Eq.~(\ref{0v(8.10)}).

It is interesting to observe that Eq.~(\ref{0v(8.14)}) indicates that in order to obtain the optimal average thickness $h_\alpha$, the corresponding variable thickness $H_\alpha({\bf y})$  has been averaged with the weight function $\rho_*({\bf y})$ given by (\ref{0v(8.11b)}).

Finally, note that in the case of compressible layers, the optimal value of the average thickness $h_\alpha$ coincides with the simple average of $H_\alpha({\bf y})$. 

\section{Discussion}
\label{0vSectionD}

Let us compare the optimal condition (\ref{0v(8.15H)}) with the condition derived in \cite{Argatov2011mb} for the problem of gap variation, for which the limit problem for the variation of the contact pressure density takes the form
\begin{equation}
m^{-1}\Delta_y \tilde{p}({\bf y})=\tilde\varphi({\bf y}),\quad {\bf y}\in\omega, \quad
\tilde{p}({\bf y})=0,\quad {\bf y}\in\Gamma,
\label{0v(8.D1)}
\end{equation}
where $\tilde\varphi({\bf y})$ is the ellipsoidal gap variation.

In this case, the optimization condition (\ref{0v(8.6)}), which is conveniently rewritten as
$$
\iint\limits_\omega \tilde{p}({\bf y})\Delta_y \theta({\bf y})\,d{\bf y}=0,
$$
after applying the second Green's formula can be replaced with the following one:
\begin{equation}
\iint\limits_\omega \tilde\varphi({\bf y})\theta({\bf y}) \,d{\bf y}=0.
\label{0v(8.D2)}
\end{equation}

Finally, with the reference to the characteristic contact region $\omega_*$, Eq.~(\ref{0v(8.D2)}) takes the form
\begin{equation}
\iint\limits_{\omega_*} \tilde\varphi({\bf y})\theta_*({\bf y}) \,d{\bf y}=0,
\label{0v(8.D3)}
\end{equation}
where we introduced the notation
\begin{equation}
\theta_*({\bf y})=1-\frac{y_1^2}{{a_1^*}^2}-\frac{y_2^2}{{a_2^*}^2}.
\label{0v(8.D4)}
\end{equation}

Thus, comparing Eqs.~(\ref{0v(8.15H)}) and (\ref{0v(8.D3)}), we see that they differ only by their weight functions $\rho_*({\bf y})$ and $\theta_*({\bf y})$, respectively, given by (\ref{0v(8.11b)}) and (\ref{0v(8.D4)}). It is interesting to observe that the maximum of the function $\rho_*({\bf y})$ does not coincide with the center of the domain $\omega_*$, while the function $\theta_*({\bf y})$ makes an emphasis namely on the central part of $\omega_*$. 

It should be noted that the least squares optimization criterion (which is analogous to the integral (\ref{0v(8.12)}) without the weight function) derived in \cite{Argatov2011mb} from an ill-posed problem underestimates the contribution of the central part of the contact domain $\omega$. Thus, to minimize the effect of the peripheral part of $\omega$, it is required to decrease the characteristic contact region $\omega_*$.


\section{Conclusion}
\label{0vSectionC}

The present study results in asymptotic solutions to the three-dimensional unilateral contact problem for a bonded thin elastic compressible or incompressible layer of variable thickness. As the main result of the present paper, the four term asymptotic expansion (\ref{0v(7.1)}) is obtained for the normal displacement at the contact surface. In the compressible case, the derived asymptotic expansion (\ref{0v(7.4)}) for the contact pressure state recovers the corresponding solution obtained in \cite{VorovichPeninin1971} for the 2D case. 

The objective of this study was to apply an asymptotic modeling approach for evaluating the sensitivity of the asymptotic model of tibio-femoral contact due to small variations in the thicknesses of the contacting articular cartilage layers. It was found that to minimize the influence of the cartilage thickness non-uniformity on the force-displacement relationship, the effective geometrical characteristics $h_1$ and $h_2$ of articular layers, which enter the asymptotic model (\ref{0v(1.1)}), (\ref{0v(1.2)}), should be determined from the introduced optimization criterion (\ref{0v(8.12)}). 

\section*{Acknowledgements}
The financial support from the European Union Seventh Framework Programme through the Marie Curie Fellowship PIIF-GA-2009-253055 is gratefully acknowledged.
Asymptotic models of articular contact were presented at the Lappeenranta University of Technology. The author is thankful to Professor A.~Mikkola for valuable and stimulating discussion.

\end{document}